\documentclass[review,11pt]{elsarticle}
\usepackage{amsmath}
\usepackage{graphicx}
\usepackage{subcaption}

\usepackage{amsthm}
\newtheorem{theorem}{Theorem}[section]

\newtheorem{lemma}[theorem]{Lemma}
\theoremstyle{definition}

\usepackage{amssymb}

\usepackage{tcolorbox}
\usepackage{rotating}
\usepackage{multirow}

\newcommand{\wideunderline}[2][2em]{%
  \underline{\makebox[\ifdim\width>#1\width\else#1\fi]{#2}}%
}
\usepackage{longtable}
\usepackage{tabularx}

\usepackage{adjustbox}

\usepackage{makecell}
\usepackage{bbold}
\usepackage{longtable}

\usepackage{lineno,hyperref}
\modulolinenumbers[5]

\makeatletter
\def\ps@pprintTitle{%
   \let\@oddhead\@empty
   \let\@evenhead\@empty
   \def\@oddfoot{\reset@font\hfil\thepage\hfil}
   \let\@evenfoot\@oddfoot
}
\makeatother

\makeatletter
\let\ams@underbrace=\underbrace
\def\underbrace#1_#2{%
  \setbox0=\hbox{$\displaystyle#1$}%
  \ams@underbrace{#1}_{\parbox[t]{\the\wd0}{#2}}%
  }
\makeatother

\usepackage[left=3cm, right=3cm, top=2cm]{geometry}
\usepackage[linesnumbered,lined,boxed,commentsnumbered]{algorithm2e}

\biboptions{authoryear}

\begin{document}

\begin{frontmatter}

\title{A mathematical programming-based solution method for the nonstationary inventory problem under correlated demand}

\author[address1]{Mengyuan Xiang\corref{mycorrespondingauthor}}
\ead{mengyuan.xiang@ed.ac.uk}

\author[address1]{Roberto Rossi}
\ead{roberto.rossi@ed.ac.uk}

\author[address1]{Belen Martin-Barragan}
\ead{Belen.Martin@ed.ac.uk}

\author[address2]{S. Armagan Tarim}
\ead{armagan.tarim@ucc.ie}

\cortext[mycorrespondingauthor]{Corresponding author. +44(0)131 651 3226}

\address[address1]{Business School, University of Edinburgh, Edinburgh, United Kingdom}
\address[address2]{Cork University Business School, University College Cork, Ireland}

\begin{abstract}
This paper extends the single-item single-stocking location non-stationary stochastic inventory problem to relax the assumption of independent demand. We present a mathematical programming-based solution method that relaxes the assumption of demand independence between time periods in an existing piecewise linear approximation strategy under the receding horizon control framework. Our method can be solved via off-the-shelf mixed-integer linear programming solvers. It can tackle demand under various assumptions: the multivariate normal distribution, a collection of time-series processes, and the Martingale Model of Forecast Evolution. We compare against solutions via stochastic dynamic programming to demonstrate that our method leads to near-optimal solutions. 
\end{abstract}


\end{frontmatter}


\section{Introduction}
Demand process characteristic plays an important role in determining the optimal inventory policy in an inventory model \citep{sethiandcheng1997}. \cite{nealeandwillems2009} argued that a nonstationary stochastic demand is the rule, rather than the exception. Many early work focused on modelling nonstationary stochastic inventory problems under the assumption that each period's demand is independent to the demand in other periods. However, as discussed in \cite{songandzipkin1993}, environmental factors, such as economic conditions, market conditions, and any exogenous conditions, have significant effects on the demand for a product, the supply, and the cost structure. In this respect, this paper aims to {\em relax the assumption of demand independence between time periods.} 

The literature on correlated demand can be roughly categorised into two groups. The first group assumes that the demand is driven by exogenous random factors, such as economic conditions and market conditions. The Markovian model is commonly adopted \citep{songandzipkin1993,sethiandcheng1997,beyerandsethi1997}, when the decision-makers fully observe demand states, and true demand distributions are known. The second group assumes demand distributions are updated with internal information, such as past demands; to capture this,  time-series processes are adopted \citep{johnsonandthompson1975,graves1999}. This study falls into the second group, in which demand forecasts are updated based on past demands. Recently, the Martingale Model of Forecast Evolution (MMFE) has been extensively investigated, as it captures a wide range of forecasting techniques, such as time-series processes and judgement-based methods \citep{NorouziandUzsoy2014, Albeyetal2015}. The model proposed in this study can be adapted to tackle demand updates following the MMFE.

We consider the periodic-review single-item single-stocking location nonstationary inventory problem under fixed and unit ordering costs, holding costs, and penalty costs. We assume that demands in different periods follow a multivariate normal distribution and that demands in successive periods are correlated. After observing demand realisations, demand forecasts for the remaining periods are updated. Replenishment plans that use these forecasts take into account the evolution of these updates over time. 

This study presents a novel approach for tackling the nonstationary lot-sizing problem under correlated demand. We leverage two key building blocks: modelling techniques originally discussed in \cite{rossietal2015}, and results from multivariate stochastic process analysis. \cite{rossietal2015} presented a mixed-integer linear programming (MILP) model for approximating the optimal $(R, S)$ policies under the assumption of independently distributed demand. The $(R, S)$ policy fixes the timing of inventory reviews ($R$) and associated order-up-to-levels ($S$) at the beginning of the planning horizon; actual order quantities are decided upon only at the beginning of each inventory review period. 
We extend the existing modelling technique of \cite{rossietal2015}, combining it with results from multivariate stochastic process analysis. The resulting model can tackle nonstationary stochastic inventory problems under demand following a multivariate normal distribution. The model developed is operationalized under a receding horizon control framework, in which only the imminent replenishment plan is implemented, and a re-planning is done at the beginning of each period for the rest of the planning horizon. Besides, we show that our approach can be adapted to tackle demand following a collection of time-series process and the MMFE.

Our approach differs from existing research on two counts. It is the first approach that combines exiting modelling techniques with results from multivariate stochastic process analysis for tackling nonstationary stochastic inventory problems under correlated demand in the receding horizon control context. It comes with the advantage of being able to tackle demands following the multivariate normal distribution, a collection of time-series processes, and the MMFE.

Our contributions to the literature on stochastic lot-sizing are the following.
\begin{itemize}
\item We present the first mathematical programming-based solution method for tackling the nonstationary stochastic inventory problem under correlated demand.
\item We combine an existing piecewise linear approximation strategy with results from multivariate normal probability theory in the context of a receding horizon control framework to model and solve the problem via a mixed integer linear programming model that can be solved by using off-the-shelf software.
\item Our approach can tackle correlated demand under various assumptions: the multivariate normal distribution, a collection of time-series processes, and the MMFE.
\item We compare our approach against solutions obtained via stochastic dynamic programming, and demonstrate that our approach leads to near-optimal solutions, while being able to tackle larger instances.
\end{itemize}

The rest of this paper is organised as follows. Section \ref{literature} surveys literature on stochastic inventory problems under demand following Markovian processes, time-series processes, and the MMFE. Section \ref{sdp} discusses multivariate normal probability theories and the stochastic dynamic programming (SDP) formulation under correlated demand. Section \ref{sp} illustrates a stochastic optimization model of the $(R, S)$ policy. Section \ref{milp} introduces our solution method for tackling the nonstationary stochastic lot-sizing problem under correlated demand. Section \ref{computationalexperiments} presents our computational study. Finally, we draw conclusions in Section \ref{conclusion}.

\section{Literature review}\label{literature}

In his landmark study, \cite{Scarf1960} proved the optimality of $(s, S)$ policy for the stochastic inventory problem with independent demand. Many researchers attempted to relax the assumption of independence by considering Markovian demand. \cite{sethiandcheng1997} considered the discrete-time finite-horizon stochastic inventory problem with Markovian demand. They showed that a state-dependent $(s, S)$ policy is optimal in the context of total cost minimisation, which consists of fixed and proportional ordering costs, proportional holding costs, and proportional backorder costs. Under the same cost structure, \cite{beyerandsethi1997}  proved the optimality of $(s, S)$ policy from the viewpoint of long-run average cost minimisation. \cite{beyeretal1998} generalised discussions in \citep{sethiandcheng1997,beyerandsethi1997}, and proved that the $(s, S)$-type policy is optimal for the finite-horizon problem, discounted-cost infinite-horizon problem, and long-run average-cost problem. \cite{chengandsethi1999} further extended the Markovian demand models to incorporate the case in which the unsatisfied demand is lost rather than backlogged, and showed that the $(s, S)$ policy is optimal. In contrast to existing studies that focus on proving the optimality of time-dependent $(s, S)$ policy and computing the optimal expected costs, \cite{nasrandelshar2018} presented a computational framework  utilising a Markovian representation to evaluate the performance measures, which include the number of backorders, on-hand inventory, inventory position, and the ordering count process of the inventory system.

For the continuous-review infinite-horizon inventory problem with Markov-modulated demand,  \cite{songandzipkin1993} proved that the state-dependent base-stock policy is optimal when the order cost is linear in the order quantity, and the state-dependent $(s, S)$ policy is optimal if there is a fixed ordering cost.  An exact procedure and a modified value-iteration algorithm were presented to compute the optimal policy parameters. \cite{lianetal2009} studied a continuous-review model for items with an exponential random lifetime and a general Markovian renewal demand process. They derived an analytical expression for computing the expected time between successive orders, the expected minimised long-run average cost, and the $(s, S)$ policy parameters. \cite{huetal2016} investigated the periodic-review infinite-horizon inventory system with a Markov-modulated demand process under $(s, S)$ policy. They proposed a Maclaurin series analysis and a Pade approximation for computing the optimal policy parameters. The Markov-modulated demand process is also widely adopted in multi-echelon inventory systems, for example, \citep{chenandsong2001,huhandJanakiraman2008,muharremogluandTsitsiklis2008,chenetal2017}.

The collective insight of these studies is that the optimal policy, either the $(s, S)$ policy or the base-stock policy,  for a model is state-dependent to reflect the dynamics of the underlying demand environment. Demand states are fully observed by decision-makers, and true demand distributions are known. This assumption has been relaxed to partially observed through the past demand data \citep{treharneandsox2002,bayraktarandludkovski2010,malladietal2018,avcietal2020}. This study differs from  existing works by assuming that future demand distributions depend on past demands, and thus demand states are unknown until past demands are observed. 

Early research on demand forecast updates generally adopted time-series processes, which use internal information, such as past demands, to update future demand distributions. \cite{johnsonandthompson1975} proved the optimality of the base-stock policy for the single-item periodic-review inventory system with proportional holding and stock-out costs and zero lead time under a mixed autoregressive-moving average demand process and the condition that demands fall in a certain lower and upper bounds. \cite{ray1981} derived an analytical expression for computing the reorder level for which the demand follows the autoregressive (AR) and moving average (MA) processes. Under the same demand processes, \cite{fotopoulosetal1988} presented a straightforward method based on the basis of probability arguments to approximate the reorder point and safety stock for which the lead time is arbitrary. \cite{graves1999} computed the base-stock policy for a single inventory system where the demand follows an integrated moving average (IMA) process. These existing studies tackle stochastic inventory problems under the assumption that demand follows certain types of time-series models. However, our method generalises existing studies and can tackle demand following a collection of time-series processes.

The Martingale Model of Forecast Evolution (MMFE) has recently been extensively explored in the modelling of correlated demand as it captures a wide range of forecasting techniques such as time-series processes and judgement-based methods. The MMFE assumes that demand forecasts for a number of periods in the future evolve over time as new information becomes available in each period. This framework was proposed by  \citep{Gravesetal1986, HeathandJackson1994} under the assumption that the information available to make forecasts grows as time grows, forecast updates are mean zero and uncorrelated with past observations, and forecast updates are stationary. 

 The MMFE has been extended and implemented by a number of authors in the context of inventory control models. \cite{toktayandwein2001} analysed the production-inventory system with stationary demand and presented a closed-form expression for computing the base-stock level in the context of expected steady-state holding and backorder costs minimisation. \cite{gallegoandozer2001} showed that the state-dependent $(s, S)$ policy and base-stock policy are optimal for stationary stochastic inventory systems with and without fixed ordering costs. \cite{IidaandZipkin2006} developed a functional approximation and a simulation-based method for approximating the base-stock policy. \cite{Luetal2006} developed easy-to-compute bounds to the optimal base-stock level, which generalised and improved existing bounds in the literature. These bounds were further used to construct near-optimal policies. \cite{ChenandLee2009} showed that many commonly adopted time-series processes such as the auto-regressive AR(1) model, the integrated moving average IMA(0,1,1) model, the general auto-regressive moving average ARMA model could be interpreted as special cases of the MMFE framework. 

These existing studies generally assumed that forecasts represent the conditional mean of demand given all information available at the time the forecast was made. Recently, \cite{NorouziandUzsoy2014} considered updates of conditional covariance of demand in addition to conditional mean. They showed that the optimal base-stock level depends on the conditional covariance, and the proposed approach yields significant cost reductions and effective decisions.  \cite{Albeyetal2015} integrated the MMFE developed by \cite{NorouziandUzsoy2014} to a chance-constrained stochastic optimization model in a rolling horizon setting. Computational studies, using data from a major semiconductor manufacturer, demonstrated that considering forecast evolution in the production planning model can lead to improved performance. Additional work on the multi-echelon inventory system with MMFE was conducted by \citep{dongandlee2003, Ziarnetzkyetal2018}. 

Inventory models with MMFE generally neglected the fixed ordering cost and aimed at computing the base-stock level, with the exception of \cite{gallegoandozer2001}, who focused on computing $(s, S)$ policies. This paper takes into account the fixed ordering cost and presents a solution method based on the $(R, S)$ policy. Besides, the MMFE assumes forecast updates are stationary. Our solution method can tackle demand updates following MMFEs, as well as nonstationary stochastic distributions. 

\section{A stochastic dynamic programming formulation}\label{sdp}
We consider a nonstationary stochastic lot-sizing problem over a $T$-period planning horizon. The demand ${d}_t$ in each period $t=1, \ldots,T$ is a random variable which follows a probability density function $g_{{d}_t}(\cdot)$ and a cumulative density function $G_{{d}_t}(\cdot)$. Let $\eta_t$ represent the realisation of $d_t$, and ${F}_t=(\eta_1, \ldots, \eta_{t-1})$ denote demand realisations up to the beginning of period $t$. In contrast to most exiting literature, we assume that demands in successive periods are not independently  distributed, but correlated. Thus, $d_t$ follows the conditional probability density function $g_{d_t}(\zeta_t|F_t)$, where $\zeta_t$ represents the value of random variable $d_t$. A full list of symbols is available in \ref{symbols}.

For our model, we will use some basic results from multivariate analysis. Let ${d}$ be a $n$-variate multivariate normal random variable with mean ${\tilde{d}}$ and covariance ${\Sigma}$, abbreviated $d \sim \mathcal{MVN}({\tilde{d}}, {\Sigma})$. 
Let $d=[d_1 \  d_2]^T$ be a partitioned multivariate normal random n-vector, $d_1=[d_1 \ \dots \ d_p]^T$, $d_2=[d_{p+1} \ \ldots \ d_n]^T$, with mean $\tilde{d}=[\tilde{d}_1 \ \tilde{d}_2]^T$ and covariance  
${\Sigma}=
\begin{bmatrix}
\Sigma_{11} & \Sigma_{12}\\
\Sigma_{21} & \Sigma_{22}
\end{bmatrix}$. 
Then the conditional distribution of ${d_2}$ given ${d_1}={\eta}_1$ is normally distributed with  mean ${\tilde{d}_2}+\Sigma_{21}\Sigma_{11}^{-1}({\eta}_1-{\tilde{d}_1})$ and variance $\Sigma_{22}-\Sigma_{21}\Sigma_{11}^{-1}\Sigma_{12}$ (Theorem \ref{conditionaldistribution}).

\begin{theorem}[Conditional distribution (see fro instance \cite{billingsley2008})]\label{conditionaldistribution}
\begin{align}
&\text{E}[{d_2}|{d_1}={\eta}_1]={\tilde{d}_2}+\Sigma_{21}\Sigma_{11}^{-1}({\eta}_1-{\tilde{d}_1})\\
&\text{Var}({d_2}|{d_1}={\eta}_1)=\Sigma_{22}-\Sigma_{21}\Sigma_{11}^{-1}\Sigma_{12}.
\end{align}
\end{theorem}

Regarding the dynamics of the inventory, we assume that replenishments $Q_t$ are placed at the beginning of period $t$ given the opening inventory level $I_{t-1}$ and demand realisations ${F}_t=(\eta_1, \ldots, \eta_{t-1})$, and delivered instantaneously. Note that $I_{t-1}$ represents the inventory at the beginning of period $t$ and the inventory level at the end of period $t-1$. The ordering cost $u(\cdot)$ comprises the fixed ordering cost $K$, and the unit ordering cost $c$, as illustrated in Eq. (\ref{orderingcost}). 
\begin{align}\label{orderingcost}
&u(Q_t)=\begin{cases} K+c\cdot Q_t & Q_t >0     \\
0 & Q_t=0 \end{cases}  
\end{align}
At the end of period $t$, the linear holding cost $h$ is charged on every unit carried from one period to the next, and the linear penalty cost $b$ is occurred for each unit of unmet demand. 
  
Given the above problem description, the objective is to schedule replenishment plans to minimize the expected total cost. In what follows, we formulate this problem as a stochastic dynamic program \citep{bellman1957} comprising the following elements.
\begin{enumerate}
\item {\bf State.} The state includes the inventory level $I_{t-1}$ and demand realisations ${F}_t =(\eta_1, \ldots, \eta_{t-1})$ at the beginning of  period $t$.
\item {\bf Action.} An action means to schedule a replenishment with quantity $Q_t$ at the beginning of period $t$, where $Q_t \in [0, \infty)$.
\item {\bf Expected immediate cost.} Let $f_t(I_{t-1}, {F}_t; Q_t)$ denote the expected immediate cost comprising ordering, holding, and penalty costs in period $t$, given initial inventory level $I_{t-1}$ and demand realisations ${F}_t$. 
\begin{align}\label{ffunction}
\begin{split}
f_t(I_{t-1}, {F}_t; Q_t)&=u(Q_t)+h\int_{{d}_t} \max(I_{t-1}+Q_t-\zeta_t, 0)g_{{d}_t}(\zeta_t|{F}_t)\text{d}(\zeta_t)\\
&\qquad + b \int_{{d}_t} \max(\zeta_t-I_{t-1}-Q_t, 0)
g_{{d}_t}(\zeta_t|{F}_t)\text{d}(\zeta_t).\\
\end{split}
\end{align}
\item {\bf Objective function.} Let $C_t(I_{t-1}, {F}_t)$ denote the expected total cost of an optimal policy over periods $t, \ldots, T$ with initial inventory level $I_{t-1}$ and up-to-date demand realisations ${F}_t$. Then, $C_t(I_{t-1}, {F}_t)$ can be written as, for $t=1, \ldots, T-1$,
\begin{align}
C_t(I_{t-1}, {F}_t)&=\min_{Q_t \geq 0}\{f_t(I_{t-1}, {F}_t; Q_t)+\int_{{d}_t}C_{t+1}(I_{t-1} +Q_t-\zeta_t, {F}_{t+1})g_{{d}_t}(\zeta_t|{F}_t)\text{d}(\zeta_t)\}, \notag\\
\end{align}
where
\begin{align}
&C_T(I_{T-1}, {F}_T)=\min_{Q_T \geq 0}\{f_T(I_{T-1}, {F}_T; Q_T)\}
\end{align}
represents the boundary condition.
\end{enumerate}


\section{\texorpdfstring{A stochastic optimisation model under $(R, S)$ policy}{Lg}} \label{sp}
In this section, we reformulate the stochastic dynamic program presented in Section \ref{sdp} as a stochastic optimization model under the $(R, S)$ policy, since our solution method is built upon modelling techniques introduced in \cite{rossietal2015} which operates under the $(R, S)$ policy for the independent demand case. We then combine it with results from multivariate stochastic process analysis. The resulting stochastic optimisation model computes $(R, S)$ policies in a way that is similar to the independent demand case. This model will further be approximated as an MILP model by adopting an existing piecewise linear approximation strategy \citep{rossietal2015} in the next section, which can be easily solved by using off-the-shelf software.


Under the $(R, S)$ policy, the timing of inventory reviews $R$ and respective order-up-to-levels $S$ are fixed simultaneously at the beginning of the planning horizon, and actual ordering quantities are decided at the beginning of each inventory review period to reach order-up-to-levels. This policy provides effective means for reducing planning instability and coping with demand uncertainty \citep{bookbinderandtan1988}. Many effective methods have been proposed for the computation of $(R, S)$ policy parameters under the assumption that the demand is independently distributed,  such as \citep{bookbinderandtan1988, tarimandkingsman2004, rossietal2015}. However, this paper leverages existing modelling techniques and results from multivariate stochastic process analysis for computing $(R, S)$ policy parameters under correlated demand.

We first introduce a binary variable $\delta_t$, $t=\{1, \ldots, T\}$, which takes value $1$ if a replenishment is placed at the beginning of period $t$ and $0$ otherwise. Then, the ordering cost in Eq. (\ref{orderingcost}) can be rewritten as follows,
\begin{align}
&u(Q_t)=K\delta_t+c\cdot Q_t.
\end{align}

At the beginning of the planning horizon, the initial inventory level is $I_0$, and we observed no demand realisations, i.e., $F_t=\O$. Following the $(R, S)$ policy, the objective is to decide the timing of replenishments $\{\delta_1, \ldots, \delta_T\}$ and the corresponding order-up-to-levels $\{S_1, \ldots, S_T\}$ for the entire planning horizon so that the expected total cost ${C}_1(I_0)$ is minimized. Thus, the problem described in Section \ref{sdp} can be formulated as the following stochastic optimization model.
\begin{figure}[!h]
\begin{align}\label{objsp}
\bar{C}_1(I_0)&=\min_{
\scriptsize
\begin{array}{c}
\delta_1,\ldots,\delta_T\\
S_1, \ldots, S_T
\end{array}}
\int_{d_1}\cdots \int_{d_T} \sum_{t=1}^T\Big( K\delta_t+c\cdot Q_t+h\text{max}(I_{t-1}+Q_t-\zeta_t,0) \notag\\
&+b\text{max}(\zeta_t-I_{t-1}-Q_t,0)\Big)g_{d_1}(\zeta_1|{F}_1)\cdots g_{d_T}(\zeta_T|{F}_T)\text{d}(\zeta_1)\cdots\text{d}(\zeta_T)
\end{align}
subject to, $t=1,\ldots,T$
\begin{align}
&Q_t=(S_t-I_{t-1})\delta_t \label{sp1}\\
&I_t=I_0+\sum_{i=1}^t(Q_i-\zeta_i) \label{sp2}\\
&Q_t, S_t \geq 0, I_t \in \mathcal{R}, \delta_t\in\{0,1\} \label{sp3}
\end{align}
\caption{A stochastic optimization model of $(R, S)$ policy under correlated demand}
\label{spmodel}
\end{figure}

The objective function (\ref{objsp}) fixes the timing of replenishments and the associated order-up-to-levels once and for all at the beginning of the planning horizon. Constraints (\ref{sp1}) states that the ordering quantity $Q_t$ is equal to order-up-to-level, minus the opening inventory level at the beginning of period $t$ if an order is placed, and 0 otherwise. Constraints (\ref{sp2}) are the inventory conservation constraints. The inventory level at the end of period t  must be equal to the inventory level at the beginning of this period, plus all received replenishments, minus demand realisations up to period $t$. Constraints (\ref{sp2}) specify domains of the order quantity, order-up-to level, inventory level, and binary variable $\delta_t$.   

Note that the objective function (\ref{objsp}) in Fig. \ref{spmodel} can be simplified as 
\begin{align}\label{objsp2}
SPC_1(I_0)&=\min_{
\scriptsize
\begin{array}{c}
\delta_1,\ldots,\delta_T\\
S_1, \ldots, S_T
\end{array}
}
\int_{d_1}\cdots \int_{d_T} \sum_{t=1}^T \Big(K\delta_t+c\cdot Q_t+h\max(I_{t-1}+Q_t-\zeta_t,0) \notag\\
&+b\max(\zeta_t-I_{t-1}-Q_t,0)\Big)g_{d_1}(\zeta_1)\cdots g_{d_T}(\zeta_T)\text{d}(\zeta_1)\cdots\text{d}(\zeta_T),
\end{align}
by applying Lemma \ref{doubleexpectation}.
\begin{lemma}[{\bf Law of total expectation} \citep{weiss2006}] \label{doubleexpectation}
 If $X$ is an integral random variable (i.e. $\text{E}[X] < \infty$) and $Y$ is any random variable, not necessarily integral, on the same probability space, then
\begin{align}
\text{E}[X]=\text{E}(\text{E}(X|Y))
\end{align}
i.e., the expected value of the conditional expected value of $X$ given $Y$ is the same as the expected value of $X$. 
\end{lemma}

Existing studies, no matter they modelled the correlated demand as a time-series process or MMFE, generally update future demand distributions based on up-to-date demand realisations. However, Eq. (\ref{objsp2}) shows that the future demand distributions are updated with unconditional means and covariances. Therefore, under Lemma \ref{doubleexpectation}, the nonstationary stochastic inventory problems under correlated demand can be tackled in a way that is similar to the independent demand cases.

\section{A new approach for nonstationary inventory problems under the receding horizon control}\label{milp}
This section first formulates the stochastic optimization model presented in Section \ref{sp} as a mixed integer linear programming (MILP) model (Section \ref{milpmodel}) by using the modelling technique introduced in \citep{rossietal2014,rossietal2015}. This model is then implemented under the receding horizon control in Section \ref{rhcmodel}. The solution method presented is adapted to tackle demand updates following MMFEs and time-series processes in Section \ref{applications}.

\subsection{\texorpdfstring{An MILP model for computing {$(R, S)$} policies with correlated demand}{Lg}}\label{milpmodel}
This section formulates the stochastic optimization model in Fig. \ref{spmodel} as an MILP model for approximating $(R, S)$ policies under correlated demand, which is built upon the modelling techniques introduced in \cite{rossietal2015}. 
 
Consider a random variable $\omega$ and a scalar variable $x$. The first order loss function is defined as $L(x, \omega)=E[\max(\omega-x,0)]$, where E denotes the expected value with respect to the random variable $\omega$. The complementary first order loss function is defined as $\hat{L}(x, \omega)=E[\max(x-\omega,0)]$. Like \cite{rossietal2015}, we model non-linear holding and penalty costs by means of this function. 

Let $d_{jt}$ represent the convolution $d_j+\ldots+d_t$. Since $d_t$ follows the multivariate normal distribution $d_t \sim \mathcal{MVN}({\tilde{d}_t}, {\Sigma})$, the demand convolution $d_{jt}$ follow the normal distribution $d_{jt} \sim \mathcal{N}(\tilde{d}_{it}, \sigma^2_{jt})$, where $\tilde{d}_{it}=\tilde{d}_i+\ldots+\tilde{d}_t$ and $\sigma^2_{jt}=\mathbb{1}^{T}\sum_{jt} \mathbb{1}$, and the symbol ``$\sim$" denotes the expectation operator.

Let $\zeta_{it}$ denote the value of random variable $d_{it}$, for $t=\{i, \ldots, j\}$. Consider a replenishment cycle $i, \ldots, t$, in which the only replenishment is placed at the beginning of period $i$. The inventory level at the end of period $t$ must be equal to the order-up-to-level at the beginning of period $i$, minus the demand convolution over periods $i, \ldots, t$, i.e. $I_t=S_i-\zeta_{it}$. Therefore, the expected on-hand stock at the end of period $t$ in Eq. (\ref{ffunction}) can be reformulated in the form of the complementary of the first order loss function,
\begin{align}\label{l2}
&\hat{L}(S_i,d_{it}) = \int_{d_i}\ldots \int_{d_t}\max (S_i-\zeta_{it},0)g_{d_i}(\zeta_i)\ldots g_{d_t}(\zeta_t)\text{d}(\zeta_i)\ldots, \text{d}(\zeta_t).
\end{align} 
And, the expected back-order at the end of period $t$ can be reformulated in the form of first order loss function,
\begin{align}\label{l1}
&L(S_i,d_{it}) = \int_{d_i}\ldots \int_{d_t}\max (\zeta_{it}-S_i,0)g_{d_i}(\zeta_i)\ldots g_{d_t}(\zeta_t)\text{d}(\zeta_i)\ldots, \text{d}(\zeta_t).
\end{align}

We introduce a binary variable $P_{jt}$ which is set to one if the most recent replenishment up to period $t$ was issued in period $j$, where $j\leq t$; if no replenishment occurs before or at period $t$, then we let $P_{1t}=1$, this allows us to properly account for demand variance from the beginning of the planning horizon. We observe that if $P_{jt}=1$, the closing inventory level of period $t$ must be equal to the order-up-to-level of period $j$ minus the demand convolution over periods $j, \ldots, t$, i.e. $I_t=S_j-\zeta_{jt}$. Then, following Eq. (\ref{l2}) - (\ref{l1}), the expected on-hand stock and back-order of period $t$ can be written by means of the complementary of the first order loss function and the first order loss function, $\sum_{j=1}^t\hat{L}(S_j, d_{jt})P_{jt}$, and $\sum_{j=1}^tL(S_j,d_{jt})P_{jt}$. Additionally, since period $j$ must be the only most recent order received up to period $t$, the following constraints must be satisfied.
 \begin{align}
 &\sum_{j=1}^tP_{jt}=1,&\\
 &P_{jt}\geq\delta_j-\sum_{k=j+1}^t\delta_k, & j=1, \ldots, t.
 \end{align}

We next employ the piecewise linear approximation technique proposed in \cite{rossietal2014,rossietal2015} to approximate the expected on-hand stock $\hat{L}(S_j, d_{jt})$ and back-order $L(S_j,d_{jt})$. This technique requires first to partition the support $\Omega$ of $d_{jt}$ into $W$ disjoint subregions $\Omega_1, \ldots, \Omega_W$. We define the probability mass $p_i=$Pr$[d_{jt} \in \Omega_i]$, and the conditional expectation $\text{E}[d_{jt}|\Omega_i]$ with associated region $\Omega_i$, for $i=1, \ldots, W$. Then, the Edmundson-Madansky upper bound can be applied to the expected on-hand stock and back-order.\footnote{Similarly, the Jensen's lower bound can be applied for approximating the expected excess inventory and back-orders as well, details  refer to Rossi. et al. (2015).} Let $\tilde{H}_t\geq 0$ and $\tilde{B}_t \geq 0$ denote the upper bounds to the expected on-hand and back-order stocks at the end of period $t$, then they are formulated as follows,
\begin{align}
&\tilde{H}_t\geq \sum_{j=1}^tS_jP_{jt}\sum_{k=1}^ip_k+\sum_{j=1}^t\left(e_W^{jt}-\sum_{k=1}^ip_k\text{E}[d_{jt}|\Omega_i]\right)P_{jt},\\
&\tilde{B}_t \geq -\tilde{I}_t+\sum_{j=1}^tS_jP_{jt}\sum_{k=1}^ip_k+\sum_{j=1}^t\left(e_W^{jt}-\sum_{k=1}^ip_k\text{E}[d_{jt}|\Omega_i]\right)P_{jt},
\end{align}
$t=1, \ldots, T$, $i =1, \ldots, W$, and $e_W^{jt}$ denotes the approximation error. Note that $\sum_{j=1}^tS_jP_{jt}=\tilde{I}_t+\sum_{j=1}^t\tilde{d}_{jt}P_{jt}$. 

Therefore, the stochastic optimization model in Fig. \ref{spmodel} for tackling the nonstationary stochastic inventory problem under correlated demand can be formulated as an MILP model in Fig. \ref{milp3}.
\begin{figure}[!h]
\begin{tcolorbox}
\begin{align}\label{milp3-1}
&\min_{\delta_t}-cI_0+c\sum_{t=1}^T\tilde{d}+\sum_{t=1}^T(K\delta_t+h\tilde{H}_t+b\tilde{B}_t)+c\tilde{H}_T
\end{align}
Subject to, for $t=1, \ldots, T$, $j=1, \ldots, t$, and $i=1, \ldots, W$,
\begin{align}
&\tilde{I}_t+\tilde{d}_t-\tilde{I}_{t-1}\geq 0 \label{milp3-2}\\
&\delta_t=0 \rightarrow \tilde{I}_t+\tilde{d}_t-\tilde{I}_{t-1}=0 \label{milp3-3}\\
&\sum_{j=1}^tP_{jt}=1,\label{milp3-4}\\
&P_{jt}\geq\delta_j-\sum_{k=j+1}^T\delta_k, \label{milp3-5}\\
&\tilde{H}_t\geq (\tilde{I}_t+\sum_{j=1}^t\tilde{d}_{jt}P_{jt})\sum_{k=1}^ip_k+\sum_{j=1}^t\left(e_W^{jt}-\sum_{k=1}^ip_k\text{E}[d_{jt}|\Omega_i]\right)P_{jt},\label{milp3-6}\\
&\tilde{B}_t \geq -\tilde{I}_t+(\tilde{I}_t+\sum_{j=1}^t\tilde{d}_{jt}P_{jt})\sum_{k=1}^ip_k+\sum_{j=1}^t\left(e_W^{jt}-\sum_{k=1}^ip_k\text{E}[d_{jt}|\Omega_i]\right)P_{jt}, \label{milp3-7}\\
&\delta_t \in \{0,1\} \label{milp3-8}\\
&P_{jt} \in \{0,1\},  \label{milp3-9}
\end{align}
\end{tcolorbox}
\caption{An MILP model for computing $(R, S)$ policies with correlated demands}
\label{milp3}
\end{figure}

The objective (\ref{milp3-1}) decides the timing and order-up-to-level of replenishments so as to minimize the expected total cost comprising ordering, holding, and penalty costs given the initial inventory level $I_0$. Constraints (\ref{milp3-2}) ensure the non-negativity of replenishments. Constraints (\ref{milp3-3}) are indicator constraints \citep{belotti2016handling} capturing the reorder condition. Constraints (\ref{milp3-4}) indicate the most recent replenishment before period $t$ was issued in period $j$. Constraints (\ref{milp3-5}) uniquely define in which the most recent replenishment prior to $t$ took place. Constraints (\ref{milp3-6}) - (\ref{milp3-7}) are approximations of expected holding and penalty costs at the end of period $t$ utilising the first order loss function and its complementary function. Constraints (\ref{milp3-8}) - (\ref{milp3-9}) set binary variables. 

The MILP model in Fig. \ref{milp3} can be easily implemented and solved using off-the-shelf software. The timing of inventory reviews are obtained from $\delta_t$ and the corresponding order-up-to-levels $S_t$ are obtained from $\tilde{I}_t+\tilde{d}_t$.

For the special case in which demand follows a standard normal distribution, the piecewise linear approximation parameters $p_i$, $\text{E}[d_{jt}|\Omega_i]$, and $e_W^{jt}$ are provided in \cite{rossietal2014}. These parameters can be applied to general normal distributions by using the standardisation formula $\hat{L}(S_j, d_{jt})=\sigma_{jt}\hat{L}\left(\frac{S_j-\tilde{d}_{jt}}{\sigma_{jt}}, Z\right)$ in \cite{rossietal2014}, Lemma 7, where $Z$ is a standard normal random variable. Note that, since in our case demands are correlated, $d_{jt} \sim \mathcal{N}(\tilde{d}_{jt}, \sigma_{jt}^2)$, where $\tilde{d}_{jt}=\tilde{d}_j+\ldots+\tilde{d}_t$ and $\sigma^2_{jt}=\mathbb{1}^{T}\sum_{jt} \mathbb{1}$. 

\subsection{An MILP model under receding horizon control}\label{rhcmodel}
This section implements the MILP model in Fig. \ref{milp3} under a receding horizon control framework, which is widely used in the inventory control literature \citep{kilicandtarim2011, duraletal2019}. 

Our algorithm proceeds as follows. At the beginning of each period $k$, the inventory level $I_{k-1}$ and up-to-date demand realisations ${F}_k$ are observed, based on which demand distributions for the remainning periods of the planning horizon $k, k+1, \ldots, T$ are updated using Theorem \ref{conditionaldistribution}. Then, inventory review periods and respective order-up-to levels are obtained by solving the $T-k+1$ periods problem using the MILP model in Fig. \ref{milp3}. Next, the replenishment plan for the current period $k$ is put into action, i.e., an order is placed if necessary. At the end of this period, the demand realisation $\eta_k$ is observed, and the total cost over periods $1, 2, \ldots, k$ is calculated. This procedure repeats until the end of the planning horizon. We simulate this process until the stopping criterion of a maximum estimation error of the expected total cost with a given confidence probability is satisfied. 

The MILP model under the receding horizon control is structured as follows (Algorithm \ref{rhcalgorithem}). 

\noindent \textbf{Update of opening inventory level and demand forecasts (line 5-6).} We update the inventory level at the beginning of period $k$, $I_{k-1}$, and demand forecasts $\text{E}[d_t|F_k]$ and $\text{Cov}(d_t,d_{t+i}|F_k)$ for $t= k, \ldots, T$ and $i=0, \ldots, T-k$ by using the conditional distribution, Theorem \ref{conditionaldistribution}.

\noindent \textbf{Computation of optimal $(R, S)$ policies (line 7).} We obtain current  replenishment decisions $\delta_k$ and $Q_k$ by solving the MILP model (Fig. \ref{milp3}) over periods $k, \ldots, T$.

\noindent \textbf{Implementation of imminent replenishment decision (line 8).} The replenishment decision $\delta_k$ and $Q_k$ in period $k$ is put into action.

\noindent \textbf{Calculation of total cost (line 9).} We calculate the total cost over periods $1, \ldots, k$. 

\IncMargin{1em}
\begin{algorithm}[h]
\SetKwData{Left}{left}\SetKwData{This}{this}\SetKwData{Up}{up}
\SetKwFunction{Union}{Union}\SetKwFunction{FindCompress}{FindCompress}
\SetKwInOut{Input}{Input}\SetKwInOut{Output}{Output}

\Input{costs ($ordering cost$, $holding cost$, $penalty cost$), $expected demand$, $covariance$, $coefficient of correlation$, $initialStock$, $errorThreshold$, $confidenceProbability$}
\Output{$totalCost$}
\BlankLine
\While{$error > errorThreshold$}{
generate demand realisations $\eta=(\eta_1, \ldots, \eta_T)$\;
$expectedTotalCost=0$\;
\For{$k =1$ \KwTo $T$}{
update opening inventory level $I_{k-1}=I_{k-2}+Q_{k-1}-\eta_{k-1}$\;
update demand forecasts $\text{E}[d_k|{F}_k]$ and  $\text{Cov}(d_k, d_{k+i}|{F}_k)$, $i=0, \ldots, T-k$\;
solve the MILP model in Fig. \ref{milp3} with $I_{k-1}$, $\text{E}[d_k|{F}_k]$ and  $\text{Cov}(d_k, d_{k+i}|{F}_k)$\;
implement imminent $\delta_k$ and $Q_k$\;
calculate up-to-date costs $totalCost+=f_k(I_{k-1}, {F}_k; Q_k)$\;
}
calculate $error$\;
}
return $totalCost$
\caption{An MILP model under receding horizon control}\label{rhcalgorithem}
\end{algorithm}\DecMargin{1em}

The receding horizon has the advantage of removing nervousness, and thus unexpected replenishment plan changes, in response to demand realisations  \citep{kilicandtarim2011}. It differs from the conventional rolling horizon control, where a fixed length of time window rolls many times over the planning horizon. This rolling horizon control is more suitable for infinite or long-horizon planning problems. Besides, the length of the time window may significantly affect the optimal replenishment plans, as pointed out in \cite{bookbinderandtan1988}. The receding horizon control set the length of the time window to the number of periods of the rest planing horizon capturing the end-of-horizon effects that is of importance for products with short life cycles \citep{duraletal2019}.



\subsection{Extensions to the Martingale Model of Forecast Evolution}\label{applications}
The Martingale Model of Forecast Evolution (MMFE) proposed by \citep{Gravesetal1986, HeathandJackson1994} is a powerful framework that can capture both time-series models of prediction and judgemental forecasts. In this section, we demonstrate that our model presented in section \ref{rhcmodel} can be embedded to capture the MMFE. 

Under the MMFE framework, demand forecasts are assumed to be available for a certain number of periods in the future, which refers to the forecast horizon. Let $H$ be the forecast horizon. At the end of period $s$, demand forecasts for the following $H$ periods are generated. Let $D_{s,t}$, $s \leq t \leq s+H$, represent the demand forecast made in period $s$ for period $t$. $D_{s,s}=D_s$ denotes the demand realisation in period $s$, since the forecasts are made after the actual demand of period $s$ is revealed. The demand forecast $D_{s,t}$, $t>s+H$, is set to a constant $\mu$. Note that the MMFE requires the demand to be stationary \citep{HeathandJackson1994}. However, our model has the advantage of being able to tackle nonstationary demands.  \cite{HeathandJackson1994} developed two classes of models for the behaviour of forecast updates: the additive model and the multiplicative model. For exposition simplicity, we mainly focus on the additive model in the rest of this section. The multiplicative model is discussed in \citep{HeathandJackson1994,toktayandwein2001,IidaandZipkin2006,NorouziandUzsoy2014}.

Let $\epsilon_{s,t}=D_{s,t}-D_{s-1,t}$ be a random variable which represents the forecast update made at the end of period $s$ for period $t$. $\epsilon_{s,t}=0$ for $t>s+H$. Let  $\mathcal{\epsilon}_s=(\epsilon_{s,s}, \epsilon_{s, s+1}, \ldots, \epsilon_{s, s+H})$ be the forecast update vector generated at the end of period $s$. $\epsilon_{s,s}$ represents the final update made to period $s$, and $\epsilon_{s, s+H}$ denotes the first forecast update made to period $s+H$. \cite{HeathandJackson1994} assumed the update vectors $\mathbf{\epsilon}_s$ are independent and identically distributed multvariate normal random vectors with mean 0 and covariance matrix $\sum=\sigma_{i,j}$,  where $\sigma_{i,j}$ represents the covariance of $\epsilon_{s, s+i}$ and $\epsilon_{s, s+j}$, $i,j=0, 1, \ldots, H$.

Since the MMFE requires the demand to be stationary \citep{HeathandJackson1994, dongandlee2003,NorouziandUzsoy2014}, the unconditional mean of the demand process is a constant $\mu$, and the conditional mean for period $t$ given information available in period $s$ is 
\begin{align}
&D_{t}=\mu+\sum_{j=0}^H\epsilon_{t-H+j,t}.   
\label{conditionalmean}
\end{align}

\cite{toktayandwein2001} developed the following unconditional covariance between $D_t$ and $D_{t+i}$,
\begin{align}
&Cov(D_t, D_{t+i})=\sum_{j=0}^H\sigma_{j,i+j}, & i=0,1, \ldots, H. \label{unconditionalcovariance}
\end{align}
The unconditional covariance, Eq. (\ref{unconditionalcovariance}), ignores the available demand information that forecast updates provided. \cite{NorouziandUzsoy2014} incorporated this information and proposed the following conditional covariance between $D_t$ and $D_{t+i}$ given all information available at end of period $s$,
\begin{align}
&Cov(D_t, D_{t+i}|F_s)=\sum_{j=1}^{t-s}\sigma_{t-s-j,t+i-s-j}, &1\leq t-s\leq H-i+1. \label{conditionalcovariance}
\end{align}
When $s=t$, $D_t$ is totally revealed; all covariance matrix elements $cov(D_t, D_{t+i}|F_s)=0$.

The MMFE differs from the multivariate normal distribution because of the length of the forecast horizon. Given the information available at the beginning of a certain period, the MMFE framework updates demand forecasts for the next $H$ periods. In contrast, the multivariate normal distribution updates demand forecasts for the remaining periods of the planning horizon. Recall that our solution method consists of an MILP model (Fig. \ref{milp3}) and its implementation under the receding horizon control framework (Algorithm \ref{rhcalgorithem}). The MILP model tackles the nonstationary stochastic inventory problem with correlated demand in a way which is similar to the independent demand case, which requires the unconditional mean and covariance of demands. Given the information available at the beginning of the planning horizon, the unconditional mean of MMFE is a constant $\mu$, and its unconditional covariance for periods $1, \ldots, H+1$ is updated via Eq. (\ref{unconditionalcovariance}). Thus, our MILP model can be adapted to solve demands which follow MMFEs. Under a receding horizon control framework, given the information available at the beginning of period $t$, the conditional expected demand and conditional covariance for periods $t, \ldots, t+H$ are updated based on Eq. (\ref{conditionalmean}) and Eq. (\ref{conditionalcovariance}). Therefore, our solution method presented in Section \ref{rhcmodel} can be immediately adapted to tackle stochastic inventory problems with demand updates following MMFEs. 

Most commonly used time-series models such as the AR, MA, and ARMA processes are special cases of the MMFE \citep{dongandlee2003,IidaandZipkin2006,NorouziandUzsoy2014}. Therefore, our model presented in Section \ref{rhcmodel} can also be adapted to tackle demands which follow all these time-series processes.


\section{Computational experiments}\label{computationalexperiments}
This section presents a numerical study to investigate the cost performance of the solution method discussed in Section \ref{milp}. We assume that the demand follows the multivariate normal distribution featuring correlations over different periods. Note that our approach comes with the advantage of tackling demands following a collection of time-series processes and the MMFE built upon the multivariate normal distribution; we thus only consider demands following multivariate normal distributions in this section.  Numerical experiments are conducted using IBM ILOG CPLEX Optimization Studio 12.7 and Eclipse 4.7.3 on a 3.2GHz Intel(R) Core(TM) with 8GB of RAM, and MATLAB R2019a on a 3.6GHz Intel(R) Xeon(R) with 48GB of RAM.

We consider nine demand patterns: two life cycle patterns (LCY1 and LCY2), two sinusoidal patterns (SIN1 and SIN2), four empirical patterns (EMP1, $\ldots$, EMP4), and a stationary pattern (STA). These demands are normally distributed with means presented in Table  \ref{demandpatterns1} and coefficients of variation $c_v \in \{0.1, 0.2\}$ (note that $\sigma_{d_t}=c_v \cdot \tilde{d}_t$). The fixed ordering cost $K$, proportional ordering cost $c$ and penalty cost $b$ range in $\{150, 300\}$, $\{0, 1\}$ and $\{5, 10, 20\}$. The coefficient of correlation $\rho \in \{-0.75, -0.25, 0.25, 0.75\}$. The proportional holding cost $h=1$.

We utilise the SDP model discussed in Section \ref{sdp} as a benchmark. Note that the SDP is pseudo-polynomial; therefore, an increase in the average value of the demand or its standard deviation will lead to a dramatic increase in the state space boundaries and computational times. We therefore adopt the Latin Hypercube Sampling \citep{mckayetal2000}, which works as follows. We generate a 9-by-5 matrix containing a Latin Hypercube Sample of 9 values on each of the five variables for each demand pattern. These five variables represent the fixed ordering cost, unit ordering cost, penalty cost, coefficient of variation, and coefficient of correlation. Each variable has ten stratifications. However, we only have two levels of fixed ordering cost, unit ordering cost, and coefficient of variation. So if the generated points are less than 0.5, then we choose the first level; otherwise, the second level. A similar approach is followed by the penalty cost and coefficient of correlation. Thus, we obtain ten instances for each demand pattern. In total, we consider 90 instances in the computational study. These parameters are presented in \ref{computationalStudyParameters}. 

We first implement the SDP model in Matlab 2019a. Since the demand follows the multivariate normal distribution, we adopt a discretisation step size of $0.5$ and a continuity correction factor of $0.25$ for demand patterns LCY1, LCY1, SIN1, SIN2, and STA to guarantee the model accuracy, and a discretisation step size of $1$ and a continuity correction factor of $0.5$ for EMP demand patterns (EMP1, EMP2, EMP3, EMP4) to ensure a reasonable computational time. Given a set of demand realisations, future demand distributions are updated based on the conditional distribution (Theorem \ref{conditionaldistribution}). Since demands are correlated over different periods, it is computationally impractical to find the optimal solution using SDP \citep{levietal2008, nasrandelshar2018}. Therefore, we propose a special structure $\sigma_{ij}^2=\rho^{|j-i|}\sigma_i\sigma_j$ to construct the covariance matrix. Accordingly, the computation of the conditional mean only requires knowledge of the realised demand in the previous period, which in return reduces the size of the state space.

We then implement the MILP model (Fig. \ref{milp3}) in the IBM ILOG CPLEX Optimization Studio 12.7 to obtain $(R, S)$ policies, and simulate these policies $1$ million times to obtain the expected total costs. We further implement the MILP model under the receding horizon control framework (Algorithm \ref{rhcalgorithem}) in Eclipse 4.7.3, which updates the $(R, S)$ policies at the beginning of each time period based on up-to-date demand realisations. We use `RHC' to represent Algorithm \ref{rhcalgorithem} presented in Section \ref{rhcmodel}. Since we operate under the assumption of normality, our models can be readily linearized by using the piecewise linearization parameters available in \cite{rossietal2014}. Specifically, we employ eleven segments in the piecewise-linear approximations of the expected on-hand stocks and back-orders. For the RHC, we adopt the stop criterion of a maximum estimation error of $0.02\%$ of the expected total cost with $98\%$ confidence probability. We present the expected total cost obtained via RHC in  \ref{computationalStudyParameters}, and compare the expected total costs obtained via the MILP model and RHC against that obtained via SDP.

\begin{figure}[!htbp]
\centering
\includegraphics[scale=0.7]{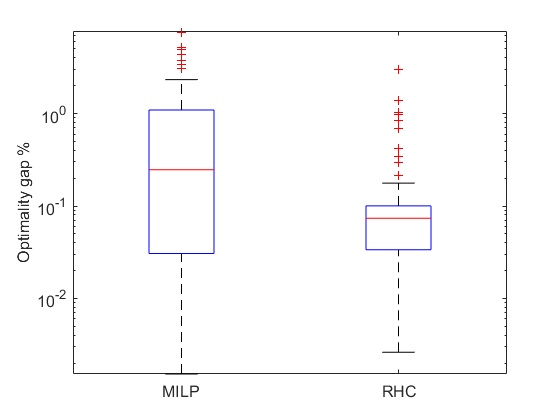}
    \caption{Optimality gaps $\%$ of our mathematical programming-based solution methods}
    \label{boxplot}
\end{figure}

We present box plots of the optimality gaps of the MILP model and RHC in Fig. \ref{boxplot}. Note that the y-axis is displayed in logarithmic scale. The average optimality gap of MILP and RHC are $0.79\%$ and $0.16\%$, and their worst-case performance are $7.45\%$ and $3.02\%$, respectively. These differences are expected due to the approximate nature of both these methods, as reported in \cite{duraletal2019}. We observe that the implementation of the receding horizon control framework can provide tighter optimality gap and reduce the dispersion of the optimality gap. 

We further present optimality gaps of the RHC for different pivoting parameters in Table \ref{optimalityGap}. Since the distribution of the optimality gap is not symmetric, and there are outliers, we report the median and interquartile range (IQR) of the optimality gap.  Regarding various demand patterns, it is difficult to draw a general remark. As we move from LCY1 to LCY2, in which demand patterns move from stable to less stable, the median slightly goes up by $0.1\%$, while the IQR goes down from $0.08\%$ to $0.03\%$. When the demand pattern changes from SIN1 to SIN2, in which demand pattern becomes wildly fluctuating, the median decreases from $0.07\%$ to $0.05\%$, while the IQR remains unchanged. 

\begin{table}[!h]
  \centering
    \begin{tabular}{|l|l|r|r|}
    \hline
    \multicolumn{2}{|l|}{Settings} & \multicolumn{1}{c|}{Median} & \multicolumn{1}{c|}{IQR} \\
    \hline
    \multicolumn{1}{|l|}{\multirow{9}[2]{*}{demand patterns}} & LCY1  & 0.08  & 0.08 \\
          & LCY2  & 0.09  & 0.03 \\
          & SIN1  & 0.07  & 0.05 \\
          & SIN2  & 0.05  & 0.05 \\
          & STA   & 0.09  & 0.07 \\
          & EMP1  & 0.08  & 0.07 \\
          & EMP2  & 0.09  & 0.20 \\
          & EMP3  & 0.05  & 0.04 \\
          & EMP4  & 0.07  & 0.08 \\
    \hline
    \multicolumn{1}{|l|}{\multirow{2}[2]{*}{fc}} & 150   & 0.08  & 0.07 \\
          & 300   & 0.07  & 0.06 \\
    \hline
    \multirow{2}[2]{*}{uc} & 0     & 0.08  & 0.03 \\
          & 1     & 0.07  & 0.07 \\
    \hline
    \multirow{3}[2]{*}{pc} & 5     & 0.08  & 0.06 \\
          & 10    & 0.07  & 0.07 \\
          & 20    & 0.08  & 0.05 \\
    \hline
    \multicolumn{1}{|l|}{\multirow{2}[2]{*}{cv}} & 0.1   & 0.07  & 0.06 \\
          & 0.2   & 0.08  & 0.07 \\
    \hline
    \multirow{4}[2]{*}{$\rho$} & -0.75 & 0.08  & 0.05 \\
          & -0.25 & 0.07  & 0.06 \\
          & 0.25  & 0.07  & 0.09 \\
          & 0.75  & 0.08  & 0.11 \\
    \hline
    \multicolumn{2}{|l|}{Overall} & 0.07  & 0.01 \\
    \hline
    \end{tabular}%
  \caption{Optimality gaps \% for different pivot parameters}
  \label{optimalityGap}%
\end{table}%

When the fixed ordering cost increases from $150$ to $300$ , the median and IQR of the optimality gap decreases from $0.08\%$ to $0.07\%$ and $0.07\%$ to $0.06\%$, respectively. This is probably due to a higher fixed ordering cost is more likely to induce longer replenishment cycles. Thus, fewer and more stable replenishments, which reduces the costs. When the unit ordering cost rises from $0$ to $1$, the median drops from $0.08\%$ to $0.07\%$, and the IQR rises from $0.03\%$ to $0.07\%$. The penalty cost has no significant correlation with the optimality gap. For instance, as the penalty cost goes up from $5$ to $20$, the median of the optimality gap remains at $0.08\%$. 

We also observe that when the coefficient of variation increases from $0.1$ to $0.2$ and the coefficient of covariance increase from $-0.75$ to $0.75$, their medians don't change significantly, but their IQRs increase. Especially, the IQR grows from $0.05\%$ to $0.11\%$ as the coefficient of covariance changes from $-0.75$ to $0.75$. This increasing trend is probably due to a higher correlation coefficient resulting in a higher demand variation, which means less stable of replenishments. 

Overall, the RHC yields a median of  $0.07\%$ and an IQR of $0.01\%$. Computational experiments demonstrate that our model leads to near-optimal solutions and low dispersions.

\section{Conclusion}\label{conclusion}
This paper considered the single-item single-stocking location nonstationary stochastic inventory problem under correlated demand. We presented the first mathematical programming-based solution method to model and solve the problem by combining an existing piecewise linear approximation strategy with results from multivariate stochastic process analysis in the context of horizon control framework. The resulting MILP model can be easily solved by using off-the-shelf software. We also showed that our solution method can be adapted to tackle demand under various assumptions: the multivariate normal distribution, a collection of time-series processes, and the MMFE.

We conducted a numerical study including 90 instances extrapolated from a test bed in the literature by means of a stratified sampling strategy. We compared the expected total costs obtained via our solution method with that via the SDP. Computational experiments showed that our mathematical programming-based solution method yields an average optimality gap of $0.16\%$ and an IQR of $0.01\%$, which demonstrate that our approach leads to near-optimal solutions and low dispersions. 

\bibliography{elsarticle-template}

\appendix 
\section{List of symbols}\label{symbols}
In this section we present a list of symbols used in this paper.

\begin{longtable}{ p{.20\textwidth}  p{.80\textwidth} } 
\hline
$T$& periods in the planning horizon\\
$d_t$ & random variable\\
$\zeta_t$ & value of random variable $d_t$\\
$\tilde{d}_t$ & the expected value of random variable $d_t$\\
$g(\cdot)$& probability density function\\
${F}_t$ & realised demand set at the beginning of period $t$\\
$I_t$ & inventory level at the end of period $t$\\
$I_0$ &initial inventory level at the beginning of the planning horizon\\
$Q_t$ & ordering quantity placed at the beginning of period $t$\\
$u(\cdot)$ & ordering cost\\
$K$& fixed ordering cost\\
$c$& proportional ordering cost\\
$h$& proportional holding cost\\
$b$& proportional penalty cost\\
$f_t(I_{t-1}, {F}_t; Q_t)$ & {immediate cost of period $t$ with opening inventory level $I_{t-1}$, realised demand set ${F}_t$, and order quantity $Q_t$}\\
$C_t(I_{t-1},{F}_t)$& the expected total cost of an optimal policy over period $t, \ldots, T$ with opening inventory level $I_{t-1}$ and realised demand set ${F}_t$\\
$S_t$& order-up-to-level of period $t$\\
$\delta_t$& binary variable\\
$\bar{C}_1(I_0)$ &expected total cost over period $1, \ldots, T$ under $(R,S)$ policy with initial inventory level $I_{0}$\\
$d_{jt}$ & a random variable denotes the demand over period $j, \ldots, t$, i.e. $d_{jt}=d_j+\ldots+d_t$\\
$\zeta_{jt}$ & value of random variable $d_{jt}$\\
$\tilde{d}_{jt}$ & expected value of the convolution $\tilde{d}_j+\ldots+\tilde{d}_t$\\
$\omega$&a random variable\\
$x$&a scalar value\\
$L(x, \omega)$&first order loss function\\
$\hat{L}(x, \omega)$& complementary of first order loss function\\
$P_{jt}$& a binary variable which is set to one if the most recent replenishment up to period $t$ was issued in period $j$, where $j\leq t$ --- if no replenishment occurs before or at period $t$, then we let $P_{1t}=1$, this allows us to properly account for demand variance from the beginning of the planning horizon\\
$\Omega$& support of $d_{jt}$\\
$W$&number of regions in a partition of $\Omega$\\
$i$ & region index ranging in $1, \ldots, W$\\
$\Omega_i$& the $i^{th}$ subregion of $\Omega$\\
$p_i$& $Pr(d_{jt} \in \Omega_i)$\\
$\text{E}[d_{jt}|\Omega_i]$ & conditional expectation of $d_{jt}$ in $\Omega_i$\\
$\tilde{H}_t$& the upper bound to the true value of $\sum_{j=1}^t\hat{L}(S_j,d_{jt})P_{jt}$\\
$\tilde{B}_t$& the upper bound to the true value of $\sum_{j=1}^tL(S_j,d_{jt})P_{jt}$\\
$e_W^{jt}$& approximation error\\
$\sigma_{jt}$ & the standard deviation of $d_{jt}$\\
$Z$& a standard normal random variable\\
$\mathcal{N}(\mu, \sigma^2)$& a normal random variable with mean $\mu$ and variance $\sigma^2$\\
$H$ & forecast horizon\\
$D_{s,t}$ & demand forecasts made in period $s$ for period $t$\\
$\epsilon_{s,t}$ & the forecast update made at the end of period $s$ for period $t$\\
$\epsilon_s$ & forecast update vector generated at the end of period $s$\\
\hline
\caption{A list of symbols} 
\end{longtable}
 
\section{Expected demands of the computational study}\label{testbed}
\begin{longtable}{|c|ccccccccc|}
\hline
\multicolumn{1}{|c|}{time period} & LCY1  & LCY2  & SIN1  & SIN2  & EMP1  & EMP2  & EMP3  & EMP4  & STA \\
\hline
    1     & 15    & 3     & 15    & 12    & 3     & 2     & 6     & 9     & 10 \\
    2     & 16    & 6     & 4     & 7     & 8     & 12    & 7     & 3     & 10 \\
    3     & 15    & 7     & 4     & 7     & 13    & 14    & 4     & 11    & 10 \\
    4     & 14    & 11    & 10    & 10    & 22    & 25    & 6     & 11    & 10 \\
    5     & 11    & 14    & 18    & 13    & 12    & 20    & 8     & 26    & 10 \\
    6     & 7     & 15    & 4     & 7     & 8     & 13    & 16    & 27    & 10 \\
    7     & 6     & 16    & 4     & 7     & 11    & 10    & 6     & 11    & 10 \\
    8     & 3     & 15    & 10    & 12    & 5     & 16    & 24    & 11    & 10 \\
\hline
\label{demandpatterns1}
\end{longtable}

\section{Computational study results}\label{computationalStudyParameters}
We present the expected total cost of our computational study in the table below. Note that ``etc" represents the expected total cost obtained via RHC ( Algorithm \ref{rhcalgorithem}). 

\footnotesize
\begin{longtable}{|c|rrrrrr|c|rrrrrr|}
    \hline
    \multicolumn{1}{|c|}{demand pattern} & \multicolumn{1}{c}{$K$} & \multicolumn{1}{c}{$c$} & \multicolumn{1}{c}{$b$} & \multicolumn{1}{c}{$cv$} & \multicolumn{1}{c}{$\rho$} &\multicolumn{1}{c|}{etc} & \multicolumn{1}{c|}{demand pattern} & \multicolumn{1}{c}{$K$} & \multicolumn{1}{c}{$c$} & \multicolumn{1}{c}{$b$} & \multicolumn{1}{c}{$cv$} & \multicolumn{1}{c}{$\rho$} &\multicolumn{1}{c|}{etc}\\
    \hline
    \endfirsthead
             
    \multicolumn{14}{c}{{\bfseries \tablename\ \thetable{} -- continued from previous page}} \\
    \hline
    \multicolumn{1}{|c|}{demand pattern} & \multicolumn{1}{c}{$K$} & \multicolumn{1}{c}{$c$} & \multicolumn{1}{c}{$b$} & \multicolumn{1}{c}{$cv$} & \multicolumn{1}{c}{$\rho$} &\multicolumn{1}{c|}{etc} & \multicolumn{1}{c|}{demand pattern} & \multicolumn{1}{c}{$K$} & \multicolumn{1}{c}{$c$} & \multicolumn{1}{c}{$b$} & \multicolumn{1}{c}{$cv$} & \multicolumn{1}{c}{$\rho$} &\multicolumn{1}{c|}{etc}\\
    \hline 
    \endhead
        
    \hline \multicolumn{14}{|r|}{{Continued on next page}} \\ 
    \hline
    \endfoot
        
    \hline 
    \endlastfoot
        
    \multirow{10}[2]{*}{LCY1} & 150   & 0     & 20    & 0.1   & -0.25 & 396.71 & \multirow{10}[2]{*}{LCY2} & 150   & 1     & 5     & 0.2   & -0.75 & 508.61 \\*
        & 300   & 0     & 5     & 0.1   & 0.75  & 539.71 &       & 300   & 0     & 20    & 0.1   & 0.25  & 693.21 \\
        & 300   & 1     & 10    & 0.2   & 0.75  & 713.81 &       & 300   & 0     & 10    & 0.2   & -0.25 & 653.81 \\
        & 150   & 1     & 5     & 0.2   & -0.75 & 460.01 &       & 300   & 1     & 20    & 0.1   & 0.25  & 781.91 \\
        & 300   & 1     & 5     & 0.1   & -0.75 & 606.21 &       & 150   & 1     & 5     & 0.2   & 0.25  & 519.21 \\
         & 150   & 1     & 20    & 0.1   & 0.25  & 499.51 &       & 300   & 1     & 10    & 0.2   & -0.25 & 738.11 \\
        & 150   & 1     & 10    & 0.1   & 0.75  & 503.11 &       & 300   & 0     & 20    & 0.1   & -0.75 & 671.61 \\
        & 300   & 0     & 10    & 0.2   & -0.25 & 558.11 &       & 150   & 0     & 5     & 0.1   & 0.75  & 435.11 \\
        & 150   & 0     & 20    & 0.2   & 0.25  & 452.61 &       & 150   & 0     & 5     & 0.1   & -0.75 & 432.61 \\
        & 300   & 0     & 5     & 0.2   & -0.25 & 538.41 &       & 150   & 1     & 20    & 0.2   & 0.75  & 603.81 \\
    \hline
    \multirow{10}[1]{*}{SIN1} & 150   & 0     & 10    & 0.2   & 0.75  & 422.51 & \multirow{10}[1]{*}{SIN2} & 150   & 0     & 5     & 0.1   & -0.25 & 394.21 \\*
        & 150   & 1     & 10    & 0.2   & 0.25  & 471.61 &       & 300   & 1     & 20    & 0.2   & 0.75  & 728.11 \\
        & 300   & 0     & 20    & 0.1   & -0.75 & 543.61 &       & 300   & 1     & 20    & 0.2   & -0.75 & 663.21 \\
        & 300   & 1     & 5     & 0.2   & 0.25  & 588.61 &       & 150   & 1     & 5     & 0.2   & -0.25 & 463.61 \\
        & 300   & 1     & 20    & 0.1   & -0.75 & 613.01 &       & 150   & 1     & 10    & 0.1   & 0.25  & 498.61 \\
        & 300   & 1     & 10    & 0.1   & -0.25 & 603.81 &       & 150   & 0     & 10    & 0.2   & -0.75 & 424.41 \\
        & 150   & 1     & 5     & 0.2   & 0.75  & 453.11 &       & 300   & 0     & 5     & 0.1   & 0.25  & 546.31 \\
        & 300   & 0     & 20    & 0.1   & 0.25  & 554.81 &       & 150   & 1     & 20    & 0.1   & 0.75  & 516.21 \\
        & 150   & 0     & 20    & 0.1   & -0.25 & 398.61 &       & 300   & 0     & 10    & 0.2   & -0.25 & 580.81 \\
        & 150   & 0     & 5     & 0.2   & -0.25 & 370.51 &       & 300   & 0     & 10    & 0.1   & 0.25  & 575.91 \\
    \hline
    \multirow{10}[1]{*}{EMP1} & 300   & 1     & 20    & 0.1   & 0.25  & 1044.01 & \multirow{10}[1]{*}{EMP2} & 300   & 0     & 10    & 0.2   & 0.75  & 1035.61 \\*
        & 150   & 0     & 20    & 0.2   & 0.25  & 641.81 &       & 300   & 1     & 10    & 0.2   & -0.25 & 1233.21 \\
        & 150   & 1     & 10    & 0.2   & 0.75  & 782.21 &       & 150   & 0     & 20    & 0.2   & -0.25 & 751.11 \\
        & 300   & 1     & 10    & 0.1   & -0.75 & 937.21 &       & 150   & 1     & 5     & 0.1   & 0.75  & 855.21 \\
        & 300   & 0     & 5     & 0.1   & 0.75  & 719.91 &       & 300   & 0     & 5     & 0.2   & 0.25  & 923.21 \\
        & 150   & 1     & 10    & 0.2   & -0.25 & 755.41 &       & 300   & 1     & 20    & 0.1   & -0.75 & 1212.51 \\
        & 300   & 0     & 5     & 0.1   & -0.75 & 701.11 &       & 150   & 0     & 5     & 0.1   & -0.75 & 623.31 \\
        & 300   & 1     & 10    & 0.2   & -0.75 & 949.81 &       & 300   & 1     & 5     & 0.2   & 0.25  & 1126.31 \\
        & 150   & 0     & 5     & 0.2   & -0.25 & 549.01 &       & 150   & 1     & 20    & 0.1   & 0.25  & 945.01 \\
         & 150   & 0     & 20    & 0.1   & 0.25  & 584.71 &       & 150   & 0     & 10    & 0.1   & -0.25 & 669.41 \\
    \hline    
    \multirow{10}[2]{*}{EMP3} & 150   & 1     & 10    & 0.2   & -0.75 & 733.71 & \multirow{10}[2]{*}{EMP4} & 150   & 1     & 5     & 0.1   & 0.75  & 821.91 \\*
        & 150   & 0     & 20    & 0.1   & 0.25  & 583.51 &       & 300   & 1     & 10    & 0.1   & -0.25 & 1135.31 \\
        & 300   & 0     & 5     & 0.1   & -0.25 & 848.41 &       & 150   & 0     & 10    & 0.2   & 0.25  & 685.51 \\
        & 300   & 0     & 5     & 0.1   & 0.75  & 854.11 &       & 300   & 0     & 20    & 0.2   & 0.25  & 1026.11 \\
        & 300   & 0     & 10    & 0.2   & 0.25  & 906.21 &       & 150   & 0     & 5     & 0.2   & -0.75 & 602.41 \\
        & 300   & 1     & 10    & 0.2   & -0.25 & 1047.41 &       & 300   & 1     & 10    & 0.1   & 0.75  & 1159.61 \\
        & 150   & 1     & 20    & 0.2   & -0.75 & 754.71 &       & 300   & 0     & 20    & 0.1   & -0.75 & 910.31 \\
        & 300   & 0     & 10    & 0.2   & 0.25  & 906.21 &       & 150   & 1     & 20    & 0.2   & -0.25 & 905.61 \\
        & 150   & 1     & 20    & 0.1   & 0.75  & 745.71 &       & 300   & 0     & 5     & 0.2   & -0.75 & 901.41 \\
        & 150   & 1     & 5     & 0.1   & -0.25 & 696.41 &       & 150   & 1     & 20    & 0.1   & 0.25  & 873.11 \\
    \hline
    \multirow{5}[2]{*}{STA} & 300   & 0     & 20    & 0.2   & -0.25 & 617.01 & \multirow{5}[2]{*}{STA} & 150   & 1     & 20    & 0.1   & -0.25 & 519.61 \\*
        & 150   & 1     & 5     & 0.1   & 0.25  & 476.51 &       & 150   & 1     & 5     & 0.2   & -0.25 & 479.31 \\
        & 300   & 1     & 10    & 0.1   & -0.75 & 663.81 &       & 300   & 0     & 5     & 0.2   & 0.75  & 572.91 \\
        & 300   & 1     & 10    & 0.2   & 0.75  & 719.21 &       & 300   & 0     & 20    & 0.2   & 0.25  & 638.71 \\
        & 150   & 0     & 10    & 0.1   & 0.75  & 443.71 &       & 150   & 0     & 5     & 0.1   & -0.75 & 401.51 \\
\end{longtable}%

\end{document}